\documentclass[8pt,openany]{article} 
\usepackage{amssymb} 
\usepackage{amsmath,amscd,mathrsfs}
\usepackage{amsmath}
\makeatletter
\renewcommand*\env@matrix[1][*\c@MaxMatrixCols c]{%
  \hskip -\arraycolsep
  \let\@ifnextchar\new@ifnextchar
  \array{#1}}
\makeatother

\usepackage{amsmath}
\usepackage{esint}
\usepackage{enumitem}
\usepackage{hyperref}
\usepackage[active]{srcltx}

\newtheorem{thm}{Theorem}

 \newtheorem{cor}{Corollary}
\newtheorem{lem}{Lemma} \newtheorem{rem}{Remark}
\newtheorem{prop}{Proposition}

\newtheorem{quest}{Question}
     \def\R{{ \! \rm \ I\!R}}

   \def \square{\hbox {$\sqcup
    $\llap {$\sqcap $}}}

\title{A simple generalization  of Garsia's conjecture}
\author{Leonardo A. Cano
Garc\'{i}a}
\date{\small{\it Mathematics department,Universidad Nacional de Colombia,\\
    Bogot\'a, Colombia.}}
\begin{document}
\maketitle
\begin{abstract}
We propose a natural generalization of a conjecture by Garsia, originally concerning the realization of conformal classes of genus-1 surfaces via embeddings in three-dimensional Euclidean space. This generalized conjecture is formulated within the framework of connections on principal bundles. We address this conjecture by providing solutions in several topologically rigid cases. As an outcome, our work yields novel parameterizations of the moduli space of conformal classes of compact surfaces of genus 1, each endowed with a clear geometric interpretation.
\end{abstract}

In~\cite{Garsia}, Garsia proved that every conformal class of a compact surface of genus $1$ can be realized via an embedding into $\mathbb{R}^3$. Specifically, for any Riemannian metric $g$ on $S^1 \times S^1$, there exists an embedding $\varphi: S^1 \times S^1 \to \mathbb{R}^3$ such that $g$ is conformally equivalent to $\varphi^* g_0$, where $g_0$ is the canonical metric of $\mathbb{R}^3$. This positively resolved a conjecture, which we refer to as \textbf{Garsia's conjecture} in this article, although it is in fact a theorem in present time. Pinkall~\cite{Pinkall} later reproved this result using the Hopf fibration. In this article, we present a framework that offers a new perspective on Garsia's conjecture. Leveraging Pinkall's results, we introduce, or at least identify, what we believe can be interpreted as an interesting generalization of Garsia's conjecture.
\\
\\
Let $P$ be an $S^1$-principal bundle over a Riemannian manifold $M$ with a connection $1$-form $\varphi$. Let $c$ be a smooth Jordan curve of $M$. Using $\varphi$ and the metric of $M$, we can define a complex structure on the total space of the pullback bundle $c^*P$ (see Lemma~\ref{prop: torus associated to connection}). The results of~\cite{Pinkall} can be reinterpreted as follows: Let $\pi: S^3 \to S^2$ be the Hopf fibration endowed with the connection $\varphi$ inherited from the quaternions $Q$. Then, for every Riemannian metric $g$ on $S^1 \times S^1$, there exists a closed curve $c$ such that its Hopf torus (defined by $\pi^{-1}(c)$), equipped with the Riemannian metric inherited from $S^3$, is conformally equivalent to $g$. The following natural question emerges, which we interpret as a \textbf{generalization of Garsia's conjecture}:
\begin{quest}\label{Question}
Let $P$ be an $S^1$-principal bundle over a Riemannian manifold $M$ with a connection $1$-form $\varphi$. For every Riemannian metric $g$ on $S^1 \times S^1$, does there exist a smooth Jordan curve $c$ such that $g$ is conformally equivalent to the complex structure (equivalently, conformal class) induced on the total space of $c^*P$ by $\varphi$?
\end{quest}
Throughout the article, we will commit a slight abuse of notation, as $c^*P$ will denote both the total space of the pullback bundle and the pullback bundle itself.
\\
\\
Rather than aiming to fully resolve this generalized conjecture across all topological settings, the primary goal of this article is to introduce this connection-theoretic framework and demonstrate its utility through a set of fundamental, topologically rigid cases. There is an intentional contrast here between the sophisticated geometric language employed—namely, connections on principal $S^1$-bundles and the moduli space of elliptic curves—and the classical simplicity of the manifolds we analyze (such as the Euclidean plane and the sphere). This deliberate choice allows us to isolate the geometric core of Garsia’s conjecture from the topological complications that arise when $H_1(M, \mathbb{Z}) \neq 0$. As a byproduct of this initial rigid analysis, we discover novel parametrizations of the moduli space of conformal classes of genus-1 surfaces, where the parameters carry a clear and concrete geometric interpretation. We elaborate on this claim in Section~\ref{Sec: Conclusions and P}.
\\
\\
Our work can be contextualized within the broader study of elliptic curves, a field renowned for fostering interactions among diverse branches of knowledge such as algebraic geometry, differential geometry, mathematical physics, and number theory. While our approach relies on relatively classical tools of differential geometry, we believe it offers a fresh, geometrically intuitive perspective on these well-established objects. By presenting these ideas in their simplest, most rigid forms, we hope to stimulate further discussion on how connection theory can illuminate the moduli space of conformal tori. Indeed, this article represents the foundational stage of a larger project dedicated to studying the moduli space of elliptic curves from a differential-geometric viewpoint (as opposed to more traditional algebraic-geometric perspectives), a research line we also explore in~\cite{CanoStandard}.
\\
\\
The paper is organized as follows. In Section~\ref{Sec: basic DG}, we define the basic geometric objects and establish the connection between holonomy and the moduli space of conformal tori. A central theme of our treatment is the consideration of manifolds with vanishing first homology, such as $\mathbb{R}^2$ and $S^2$. We show in Section~\ref{Sec: Connection Complex Structures} that in this topological setting, the holonomy of a Jordan curve is a rigid quantity determined entirely by the curvature flux, which greatly simplifies the reconstruction of the conformal class. (We defer the treatment of cases with non-trivial homology, where this topological rigidity is lost, to a forthcoming article). In Sections~\ref{Sec:Euclidean case} and \ref{Sec: examples sphere}, we apply this rigid framework to the Euclidean plane and to line bundles over the sphere. Finally, in Section~\ref{Sec: Conclusions and P}, we discuss how this perspective provides a clear blueprint for future generalizations to surfaces with non-trivial topology.
\section{Some basic differential geometry}\label{Sec: basic DG}
In this article, we work exclusively within the smooth category. To fix terminology and make the article more accessible, we present a brief synthesis of the notion of connection on $S^1$--principal bundles in this section. The theory discussed here can be found extensively in \cite{Kobayashi1} and \cite{Kobayashi2}, and in a more synthetic way in \cite[Section 1.1]{BerlineGetzlerV}.\\
\\
We define a \textbf{torus} as any manifold diffeomorphic to $S^1 \times S^1$, and a \textbf{complex torus} as a torus equipped with a complex structure. Throughout this article, we denote by $\mathscr{M}$ the \textbf{moduli space of conformal classes of compact surfaces of genus $1$}, which is well known to be $\mathbb{H}/SL(2,\mathbb{Z})$.
\\
\\
On a $S^1$--principal bundle $\pi: P \to M$ the {\bf vertical bundle} is the subbundle $V:={\rm Ker } d \pi$ of $TP$. A {\bf connection} for $\pi: P \to M$ is the invariant choice of a subbundle $H$ such that $TP=H\oplus V$ such bundle is named {\bf horizontal bundle}.
\\
\\   
Every $S^1$--principal bundle $P$ over a manifold comes equipped with a vertical vector field $X_\theta=\partial_\theta$ defined by
$$X_\theta(p)=\partial_\theta (e^{i\theta} p)\vert_{\theta=0}.$$ 
Every connection of a $S^1$--principal bundle $P$ is described by a $1$--form $\varphi \in \Omega^1(P)$ which satisfies
\begin{itemize}
\item[i)] {\bf Verticality Condition} $\varphi (\partial_\theta)=1$.
\item[ii)]{\bf Equivariance Condition} for all $g \in S^1$,  
$$R_g^*(\varphi)_p (X)= \varphi_{pg}(dR_g(X)),$$
where $p\in P$, $X\in T_pP$, and $R_g$ is the right action of $S^1$ which we  denote  $R_g(p)= pg$.
\end{itemize} 
$1$--forms satisfying these conditions are called {\bf connection $1$--forms}.
\\
\\
A connection $1$--form $\varphi$   induces the horizontal bundle (equivalently the connection)  
$$
H_{p}=\{v-\varphi(v)\partial_\theta: v \in T_pP \}.
$$
In this framework  $d\varphi$ defines a {\it horizontal $2$--form} and hence it induces a $2$--form $\Omega$ that is {\bf the curvature of the connection $\varphi$}.
\\
\\
Let $\pi: P \to M$ be a $S^1$--principal bundle over $M$ and suppose that $\varphi$ is connection $1$--form of $P$.   Given a smooth curve $c:[a,b] \to M$, a curve $\tilde{c}:[a,b] \to P$ is a {\bf $\varphi$--horizontal lift of $c$} if $\pi \circ \tilde{c}=c$ and $\tilde{c}$ is a solution of the differential equation $\varphi(\frac{d}{dt} \tilde{c})=0$.  
\\
\\
In Section~\ref{Sec: examples sphere}, we will use the equivalence between the category of principal $S^1$-bundles and Hermitian line bundles. Let us quickly recall this equivalence. Given a line bundle $L$ with a Hermitian metric $h$, we define $$P_L:=\{v\in L: h(v,v)=1\}.$$ The right action of $S^1$ is given by $v\cdot \alpha = \alpha \cdot v$, where the latter denotes scalar multiplication by $\alpha \in S^1$ on the vector $v$ within the fiber $L_{\pi v}$. This structure $(\pi: P_L \to M)$ forms a principal $S^1$--bundle, and the map $L \mapsto P_L$ defines a functor from the category of line bundles over $M$ to the category of $S^1$ principal bundles on $M$.
\\
\\
Conversely, for a principal $S^1$-bundle $P$ over $M$, its associated complex line bundle $L_P$ is constructed as the quotient of $P \times \mathbb{C}$ by the $S^1$-action $(p, v) \cdot z = (p \cdot z, z^{-1} v)$ for $p \in P$, $v \in \mathbb{C}$, and $z \in S^1$. The projection map for $L_P$ is induced by the projection of $P$ onto $M$, and $L_P$ naturally inherits a Hermitian metric.
\\
\\
Given a principal $S^1$-bundle $\pi_P: P \to M$ and the associated complex line bundle is $L_P := (P \times \mathbb{C})/S^1$, sections $s$ of $L_P$ correspond to $S^1$-equivariant functions $\tilde{s}: P \to \mathbb{C}$ satisfying $\tilde{s}(p \cdot z) = z^{-1} \tilde{s}(p)$. More precisely, an equivariant function $\tilde{s}$ defines the section $s \in \Gamma(L_P)$ by $s(m) = [p, \tilde{s}(p)]$, where for each $m \in M$ we choose a $p \in \pi_P^{-1}(m)$. Conversely, a section $s: M \to L_P$ defines an equivariant function $\tilde{s}$ where $\tilde{s}(p)$ is the complex number such that $(p, \tilde{s}(p)) \in s(\pi_P(p))$ (see \cite[Proposition 1.7]{BerlineGetzlerV} for further details).
\\
\\
The covariant derivative $\nabla$ on $L_P$ induced by a connection $1$-form $\varphi$ on $P$ acts on a section $s$ (represented by an equivariant function $\tilde{s}$) in the direction of a vector field $X$ on $M$ as follows: we lift $X$ to a horizontal vector field $\tilde{X}$ on $P$, and then $\nabla_X s$ is the section corresponding to the equivariant function $\tilde{X}(\tilde{s})$. Locally, if $\sigma: U \to P$ is a local section of $P$, a section of $L_P|_U$ can be written as $s(m) = [\sigma(m), f(m)]$ for some function $f: U \to \mathbb{C}$. The connection $\nabla$ has a local connection $1$-form $A = \sigma^* \varphi \in \Omega^1(U; i\mathbb{R})$, and the covariant derivative acts as $\nabla_X s = [\sigma(m), X(f)(m) + i A(X)(m) f(m)]$.
\\
\\
Conversely, let $\nabla$ be a connection on a Hermitian line bundle $L$ that is compatible with the Hermitian metric $h$, i.e., $$X h(s_1, s_2) = h(\nabla_X s_1, s_2) + h(s_1, \nabla_X s_2)$$ for any vector field $X$ and sections $s_1, s_2$ of $L$. Such a connection has a local connection $1$-form $A \in \Omega^1(U; i\mathbb{R})$ with respect to a local orthonormal frame $e$ of $L|_U$, defined by $\nabla e = i A \otimes e$. The principal $S^1$-bundle $P_L$ of unit norm vectors in $L$ has a connection $1$-form $\varphi$ which, when pulled back by a local section $\sigma: U \to P_L$ (which corresponds to a local orthonormal frame $e$ of $L|_U$), yields $\sigma^* \varphi = A$.
\\
\\
These constructions define two functors: one from the category of principal $S^1$-bundles with connection to the category of Hermitian line bundles with compatible connections, and vice versa. Furthermore, the isomorphisms establishing the equivalence of the underlying categories of bundles extend to isomorphisms of these categories with connections. Therefore, we have an equivalence of categories between principal $S^1$-bundles endowed  with connection $1$--forms  and Hermitian line bundles equipped with compatible connections. These constructions establish the  equivalence of categories we will use in Section~\ref{Sec: examples sphere}.
\section{Connections and complex structures}\label{Sec: Connection Complex Structures}
Let $\pi:P \to M$ be an $S^1$--principal bundle over a Riemannian manifold $M$, equipped with a connection $1$--form $\varphi$. Consider a smooth Jordan curve $c:[0,L] \to M$. 
\begin{lem}\label{prop: torus associated to connection}
Let $c:[0,L] \to M$ be a smooth Jordan curve parametrized by arc length. Let $P$ be an $S^1$--principal bundle over the Riemannian manifold $M$, equipped with a connection $1$--form $\varphi$. Then, the total space of the pullback bundle $c^*P$ is a torus, which is endowed with a complex structure naturally induced by the connection $\varphi$ and the Riemannian metric $g$. Moreover, this total space is biholomorphic to $\mathbb{C}/\Lambda$, where $\Lambda = \mathbb{Z} + (-\alpha+iL)\mathbb{Z}$, and $\alpha$ is the angle of holonomy of $c$.
\end{lem}
{\bf Proof:}
\\
Let $\eta:[0,L] \to P$ be a $\varphi$--horizontal lift of $c$. Since $c$ is a smooth Jordan curve, it is closed, meaning $\eta(0)$ and $\eta(L)$ belongs to the same fiber. Consequently,  there exists a unique $\alpha \in \mathbb{R}$ such that $\eta(L)=\eta(0) \cdot e^{i\alpha}$. By definition, $e^{i\alpha}$ is the holonomy of the curve $c$, and $\alpha$ is its angle of holonomy.
\\
\\
The map $c$ can be naturally extended from $[0,L]$ to an $L$--periodic smooth function on $\mathbb{R}$, which we continue to denote by $c$. Similarly, its horizontal lift $\eta$ can be extended to an $L$--periodic horizontal lift of the extended $c$ on $\mathbb{R}$, which we also continue to denote $\eta$. The extended maps satisfy
\begin{equation}\label{eq:periodicidad de Chi}
\eta(s+mL)=\eta(s) \cdot e^{i\alpha m} \text{ and } c(s+mL)=c(s).
\end{equation}
Let us onsider the map $\chi: \mathbb{R}^2 \to P$ defined by
$$
\chi(s_1, s_2) := \eta(s_2) \cdot e^{is_1}
$$
The equations in (\ref{eq:periodicidad de Chi}) imply that the map $(\varphi+is) \mapsto \chi(\varphi,s)$ is periodic with respect to the lattice $\Lambda:=2\pi \mathbb{Z}+(-\alpha+iL)\mathbb{Z}$, as the following calculation shows
$$
\chi(s_1+2\pi n-m\alpha,s+mL)=e^{i(s_1+2\pi n-m\alpha)}\eta(s+Lm)=\chi(\varphi,s).
$$
We denote $[s_1+is_2]$ as the class of the complex number $s_1+is_2$ in the quotient $\mathbb{C}/\Lambda$. It is straightforward to see that the function $[s_1+is_2] \mapsto \chi(s_1,s_2)$ is a diffeomorphism from $\mathbb{C}/\Lambda$ onto its image $T:=\chi(\mathbb{R}^2)$, which we have thus shown to be a torus. The natural complex structure of $\mathbb{C}/\Lambda$ then induces a complex structure on the torus $T$. $\square$
\\
\\
The total space of the pullback bundle $c^*P$ is diffeomorphic to the torus $T:=\pi^{-1}(c([0,L]))$. Therefore, we will often refer to this complex torus, which Lemma~\ref{prop: torus associated to connection} defines for a smooth Jordan curve $c$, simply as $c^*P$. We will implicitly omit its dependence on the connection $\varphi$ and the metric $g$ on $M$.
\begin{rem}\label{rem:Parameter length and area}
The main ingredients for defining a complex torus via Lemma~\ref{prop: torus associated to connection} are an $S^1$--principal bundle $\pi: P \to M$ equipped with a connection $1$--form $\varphi$, and a Riemannian metric on the base manifold $M$. If Question~\ref{Question} is answered positively for these ingredients, then Lemma~\ref{prop: torus associated to connection} provides a first way to geometrically interpret the parameters of the moduli space $\mathscr{M}$ as the length of $c$ in the manifold $M$ and the holonomy of the curve $c$. Moreover, if we fixed the curve $c$ and varied the Riemannian metric on $M$ or the connection $\varphi$ on $P$, we could observe how the conformal class of the torus $c^*P$ changes.
\end{rem}
Here we prove a minor generalization of \cite[Theorem 1]{SingerThorpe}. For this, we use the notion of a {\bf smooth triangle} in $M$, which is an embedding $\psi:T \to M$ of a triangle $T$ that extends to a smooth map in an open neighborhood of $T$. A similar generalization of \cite[Theorem 1]{SingerThorpe} was claimed in the proof of \cite[Proposition 1]{Pinkall}.
\begin{thm}\label{thm: holonomy in terms curvat}{\upshape (cf. \cite[Theorem 1]{SingerThorpe})}
Let $P$ be an $S^1$--principal bundle over a manifold $M$, equipped with a connection given by the $1$--form $\varphi$. Let $c:[a,b]\to M$ be a smooth closed curve such that it forms the boundary of a smooth triangle $T$ in $M$ (i.e., $c = \partial T$). Then, the holonomy of $c$ is given by
$$e^{i\int_T i^* \Omega}$$
where $i:T \to M$ is the inclusion of the triangle $T$ in $M$. Moreover, the angle of holonomy is $\alpha := \int_T i^* \Omega$.
\end{thm}
{\bf Proof:}
\\
Let $T$ have vertices $v_0,v_1,v_2$. As in \cite{SingerThorpe}, we define a {\it horizontal lift} of the inclusion $i:T \to M$. This is achieved by taking horizontal lifts of all lines connecting $v_1$ to points on the segment $v_0v_2$, ensuring these lifts start from a fixed point $w_1 \in P$ such that $\pi(w_1)=v_1$.
\\
\\
From \cite[Page 192]{SingerThorpe}, we have:
\begin{equation*}
\int\int_T i^* \Omega = \int\int_T (\pi\circ \tilde{h})^* \Omega = \int\int_T \tilde{h}^* \pi^* \Omega = \int\int_T \tilde{h}^*(d\varphi).
\end{equation*}
By Stokes' Theorem, this simplifies to:
\begin{equation}\label{eq: curvature integral and angles}
\int \int_T i^* \Omega = \int_{\partial T} \tilde{h}^* \varphi.
\end{equation}
Let $\tilde{c}$ be the horizontal lift of the curve $c$ such that $\tilde{c}(v_1)=w_1 \in P$. Let $\tilde{\beta}$ denote the lift $\tilde{h}$ restricted to the boundary $\partial T$. By construction, $\tilde{\beta}$ is a horizontal lift along the edges $v_0 v_1$ and $v_1v_2$ of the triangle $T$, thus coinciding with $\tilde{c}$ along these segments.
\\
\\
For the remaining segment $v_0v_2$, there exists a smooth function $f:v_0v_2 \to \mathbb{R}$ such that
$$
\tilde{\beta}(t)=\tilde{c}(t) \cdot e^{i f(t)},
$$
for $t\in v_0v_2$, with $f(v_2)=0$. (For more details, see \cite[Lemma 2, Section 7.1]{SingerThorpe}).
\\
\\
\cite[Lemma 2, Section 7.1]{SingerThorpe} also states that along the segment $v_0v_2$,
\begin{equation*}
\varphi\left(\frac{d}{dt}\tilde{\beta}\right) = \frac{d}{dt}f.
\end{equation*}
From Equation (\ref{eq: curvature integral and angles}), which states $\int\int_{\partial T} \tilde{h}^* \varphi = \int \int_T i^* \Omega$, we can then conclude that
$$
\int \int_T i^* \Omega = -f(v_0).
$$
Here, $-f(v_0)$ represents the holonomy angle of $w_0$ around the Jordan curve $c$. $\square$
\begin{rem}\label{rmk: holonomy}
While Theorem 1 is formulated for a triangle, it carries broader topological implications. Since any smooth Jordan curve $c$ that is homologically trivial ($[c]=0 \in H_1(M, \mathbb{Z})$) acts as the boundary of a $2$-chain $\Sigma = \sum n_i T_i$, the holonomy along $c$ is uniquely determined by the integral of the curvature over $\Sigma$. In spaces such as $\mathbb{R}^2$ and $S^2$, where every Jordan curve is null-homologous, the holonomy is thus a rigid geometric quantity depending only on the enclosed region and the curvature form of the connection. In a work in preparation, we shall investigate how the holonomy varies with the connection for curves with non-trivial homology, a study which sheds further light on the generalization of Garsia's conjecture.
\end{rem}
The following trivial corollary of Theorem~\ref{thm: holonomy in terms curvat} demonstrates that Question~\ref{Question} does not always have a positive answer, thus justifying its inquiry.
\begin{cor}\label{cor:Curvature 0}
Let $\varphi$ be a connection $1$--form on the trivial principal bundle $\pi: \mathbb{R}^2 \times S^1 \to \mathbb{R}^2$ (where $\pi(x,y,\alpha)=(x,y)$), and suppose that its curvature $\Omega$ is zero everywhere. Then, for any smooth Jordan curve $c$ in $\mathbb{R}^2$, the associated complex torus $c^*P$ is biholomorphic to $\mathbb{C}/(\mathbb{Z} + iL\mathbb{Z})$, where $L$ is the length of the curve $c$.
\end{cor}
\begin{rem}\label{rem:holonomy as an integral of curvature}
Theorem~\ref{thm: holonomy in terms curvat}, within the context of its hypotheses, offers a second possible interpretation for the parameters of the moduli space $\mathscr{M}$ (or more precisely, for the subset of complex tori obtained as $c^*P$ via Lemma~\ref{prop: torus associated to connection}). These parameters are now understood as the integral of the curvature of the connection $\varphi$ on $P$ (which is precisely the holonomy of the curve $c$) and the length of $c$.
\end{rem}
In \cite{Pinkall}, Pinkall relates the length of a Jordan curve $c$ on $S^2$ to its holonomy under the natural connection of the Hopf fibration. For the sphere, the holonomy of any closed curve $c$ is determined by $A/2$, where $A$ is the \textbf{signed area} enclosed by $c$. This illustrates the rigidity discussed in Remark~\ref{rmk: holonomy}: since every Jordan curve on $S^2$ is null-homologous, the holonomy is purely a function of the enclosed area and the curvature flux. Utilizing Lemma~\ref{prop: torus associated to connection}, Pinkall shows that the complex tori associated with these connections are of the form $\mathbb{C}/(\mathbb{Z} + \tau \mathbb{Z})$, with $\tau = A/2 + iL/2$ (\cite[Proposition 1]{Pinkall}). By applying the isoperimetric inequality, which relates the length and the enclosed area of a Jordan curve, he demonstrates that the set of such complex numbers covers a fundamental domain in $\mathbb{H}$ under the action of $SL(2,\mathbb{Z})$, thereby providing a positive answer to Question~\ref{Question} for the Hopf fibration.
\\
\\
More generally, Lemma~\ref{prop: torus associated to connection} suggests that on manifolds with a vanishing first homology group, resolving Question~\ref{Question} requires a precise relationship between the length of a smooth Jordan curve $c \subset M$ and its holonomy. Specifically, the challenge lies in identifying a fundamental domain for the $SL(2,\mathbb{Z})$ action within the set of complex numbers $z = \alpha + iL$, where $L$ is the length of $c$ and $\alpha$ is its holonomy angle. In the cases studied by Pinkall, the topological rigidity of the holonomy ensures that $\alpha = A/2$; this allows the isoperimetric inequality $L^2 \geq 4\pi A - A^2$ to define the geometric constraints on the lengths and areas of smooth Jordan curves in the sphere. Consequently, the resulting set of complex numbers $A/2 + iL/2$ is shown to cover a fundamental domain for the moduli space of conformal classes of the torus. In the following section, we will extend this principle by exploring the interplay between holonomy and length for more general connections where similar geometric relations persist.
\section{The Euclidean case}\label{Sec:Euclidean case}
In this section, we study connections associated with the trivial principal bundle $P:=\mathbb{R}^n \times S^1$ over $\mathbb{R}^n$. This is, up to isomorphism, the only principal bundle over $\mathbb{R}^n$. By the Verticality Condition and Equivariance Condition explained in Section~\ref{Sec: basic DG}, the smooth connections $\varphi$ on $P$ are the $1$--forms on $P$ given in the trivialization $(x_1,\ldots,x_n, \theta)\mapsto (x_1,\ldots,x_n, e^{i\theta})$ of $P$ by
\begin{equation*}\label{eq:connection Rn}
\varphi=d\theta+\sum_{i=1}^n a_i dx_i,
\end{equation*}
where $a_i$ are smooth functions on $\mathbb{R}^n$. If $\Omega$ is the curvature of $\varphi$, then $\Omega=d\alpha$, where $\alpha:=\sum_{i=1}^n a_idx_i$. 
\\
\\
In general, we cannot directly apply Theorem~\ref{thm: holonomy in terms curvat} to every closed curve, as not every curve in a manifold is the boundary of a region homeomorphic to a triangle. However, in $\mathbb{R}^n$, we benefit from both a global trivialization and the fact that $H_1(\mathbb{R}^n; \mathbb{Z}) = 0$. This topological property ensures that every smooth Jordan curve $c$ is the boundary of a smooth surface $\Sigma$ (for instance, a Seifert surface). By Stokes' Theorem, the holonomy along $c$ is given by the integral of the curvature $2$-form $\Omega$ over $\Sigma$; as noted in Remark~\ref{rmk: holonomy}, this value remains invariant for any choice of surface $\Sigma$ bounding $c$.
\\
\\
In terms of a global trivialization, given a smooth Jordan curve $c:[0,L] \to \mathbb{R}^n$, the equation for its horizontal lift $\tilde{c}(t)=(c(t),e^{i \theta(t)})$ is
$$
\theta'(t) = \sum_{i=1}^n a_i(c(t))c_i'(t).
$$ 
Consequently, the holonomy angle of $c$ is given by the contour integral $\int_c X \cdot ds$, where $X(p):=(a_1(p),\ldots,a_n(p))$ is the vector field representing the connection forms. This leads to the following proposition:
\begin{prop}\label{prop: torus Rn}
Let $X=(a_1,a_2,a_3)$ be a smooth vector field on $\mathbb{R}^3$, interpreted as a magnetic potential. Let $P = \mathbb{R}^3 \times S^1$ be the trivial principal bundle over $\mathbb{R}^3$, endowed with the connection $1$-form $\omega = d\theta + \sum a_i dx_i$. For any smooth closed curve $c:[0,L] \to \mathbb{R}^3$ parametrized by arc length, the complex torus $c^*P$ is biholomorphic to $\mathbb{C}/\Lambda$, where the lattice is given by $\Lambda = 2\pi\mathbb{Z} + \left( \oint_c X \cdot ds + iL \right)\mathbb{Z}$.
\end{prop}
\begin{rem}\label{rem:parametrization of moduli Rn}
On the trivial principal bundle over $\mathbb{R}^n$, the angle of holonomy can be calculated directly via the natural global trivialization. This provides a geometric interpretation of the parameters of the moduli space in terms of a contour integral of the vector field associated with the connection form. Since $H_1(\mathbb{R}^n; \mathbb{Z}) = 0$, this holonomy is an invariant quantity that can be equivalently expressed as the integral of the curvature over any surface $\Sigma$ bounding the curve $c$. In the specific case described in Proposition~\ref{prop: torus Rn}, the holonomy angle for the conformal class associated with $c^*P$ is given by the line integral of the magnetic potential, suggesting that the parameters of the moduli space admit physical interpretations related to magnetic flux and gauge invariance.
\end{rem}
\subsection{Some examples in $\R^2$}\label{Sec:Connect R2}
In this section, we specialize the previous results to the trivial principal bundle $P:=\mathbb{R}^2 \times S^1$ over $\mathbb{R}^2$. A generic connection in this case is given by
\begin{equation*}\label{eq:connection R2}
\varphi=A(x,y)dx+B(x,y) dy+d\theta.
\end{equation*}
The curvature $\Omega$ of $\varphi$ is
$$\Omega=\left(\partial_x B - \partial_y A \right) dx \wedge dy.$$
By Theorem~\ref{thm: holonomy in terms curvat} and Lemma~\ref{prop: torus associated to connection}, the complex tori obtained from the connection $\varphi$ and the Euclidean metric of $\mathbb{R}^2$ are biholomorphic to $\mathbb{C}/\Lambda$. Here, $\Lambda$ is the lattice generated by $1$ and
$$ -\int \int_T \left(\partial_x B - \partial_y A \right) dx dy+iL = -\oint_{c} X \cdot ds+i L, $$
where $X$ is the vector field on $\mathbb{R}^2$ defined by $X(x,y)=(A(x,y),B(x,y))$, and the double integral is defined over the compact closed region $T$ enclosed by the curve $c$ (which exists by the Jordan curve theorem). Therefore, to answer Question~\ref{Question} positively for the connection $1$--form $\varphi$ and the Euclidean metric of $\mathbb{R}^2$, we need to find a way to identify a fundamental domain for the $SL(2,\mathbb{Z})$ action on the upper half-plane $\mathbb{H}$ within the set of complex numbers
$$ \left\{ -\oint_c X \cdot ds + iL \mid c \text{ is a smooth Jordan curve in }\mathbb{R}^2 \right\}. $$
As the following examples will show, the answer to Question~\ref{Question} depends on the vector field $X$.
\\
\\
{\bf Example 1:}
\\
If $A$ and $B$ in (\ref{eq:connection R2}) are constant functions, then the curvature $\Omega$ is $0$. In this case, Corollary~\ref{cor:Curvature 0} shows that the complex tori are conformally equivalent to {\it a rectangular tori}. This implies that it is not possible to describe every conformal class of a genus $1$ surface as a complex torus induced by Lemma~\ref{prop: torus associated to connection} for the connection $\varphi$ defined in (\ref{eq:connection R2}) and the euclidean metric of $\R^2$.
\\
\\ 
{\bf Example 2:}
\\
Let us consider the connection $1$--form $\varphi:=-ydx+xdy+d\theta$ for the trivial bundle $P:=\mathbb{R}^2 \times S^1$ over $\mathbb{R}^2$. The curvature of this connection is $2 dx \wedge dy$. By the reasoning exposed above, the complex tori obtained via Lemma~\ref{prop: torus associated to connection} from smooth Jordan curves in $\mathbb{R}^2$ are biholomorphic to $\mathbb{C}/\Lambda$. Here, $\Lambda$ is the lattice generated by $1$ and $\frac{A}{2\pi}+i\frac{L}{2\pi}$, where $L$ and $A$ are the length and signed area, respectively, of the planar region enclosed by a smooth Jordan curve $c$.
\\
\\
The only restriction on $A$ and $L$ is given by the isoperimetric inequality in $\mathbb{R}^2$, namely $L^2 \geq 4 \pi |A|$, which is equivalent to $\left(\frac{L}{2\pi}\right)^2 \geq \frac{1}{\pi} \left|\frac{A}{2\pi}\right|$. Since $A$ represents a signed area, the set of all possible complex numbers $A+iL$ (derived from Jordan curves) is considered symmetric with respect to the line $A=0$. The collection of scaled complex numbers $\frac{A}{2\pi}+i\frac{L}{2\pi}$ that satisfies this isoperimetric inequality corresponds to a region in the upper half-plane $\mathbb{H}$ that contains the fundamental domain $\mathscr{D}:=\{z=x+iy \mid -1/2 \leq x \leq 1/2, x^2 +y^2>1\}$ for the $SL(2,\mathbb{Z})$ action on $\mathbb{H}$.
\\
\\
Hence, the fundamental domain $\mathscr{D}$ is contained in the set
\begin{equation*}
\begin{split}
&\{ \frac{A}{2\pi}+i\frac{L}{2\pi} \mid L \text{ and } A \\
 &  \text{ are the length and signed area of a smooth Jordan curve in }\mathbb{R}^2  \},
\end{split}
\end{equation*}
and every conformal class of a compact surface of genus $1$ can be realized as $c^*P$ for some Jordan curve $c$ in $\mathbb{R}^2$.
\\
\\
In this context, and more generally for any case where the base manifold has a vanishing first homology group, a form of \textbf{inverse isoperimetric problem} arises. This type of problem was central to the approach in \cite{Pinkall}, where it was established that the set of complex numbers $A+iL$, derived from the area and length of smooth Jordan curves on the sphere, covers a fundamental domain for the action of $SL(2,\mathbb{Z})$ on $\mathbb{H}$.
\\
\\
{\bf Example 3:}
\\
Let's consider the connection $1$--form $\varphi:=-y^2dx+x^2dy+d\theta$ for the trivial bundle $\mathbb{R}^2 \times S^1$ over $\mathbb{R}^2$. The curvature of this connection is $\Omega:=(2x+2y) dx \wedge dy$. We can calculate the integral of the curvature $\Omega$ over the ball $B_R(a,b)$ (a disk of radius $R$ centered at $(a,b)$):
\begin{equation*}
\begin{split}
&\int \int_{B_R(a,b)}(2x+2y) dx \wedge dy= (a+b)R^2 \cdot 2\pi = 2\pi R^2(a+b).
\end{split}
\end{equation*}
Hence, the complex torus associated with a circle of radius $R$ centered at $(a,b)$, via Lemma~\ref{prop: torus associated to connection}, is conformally equivalent to $\mathbb{C}/\Lambda$. Here, $\Lambda$ is the lattice $$2\pi\mathbb{Z}+(\frac{2\pi R^2(a+b)}{2\pi}+i\frac{2\pi R}{2\pi})\mathbb{Z} = 2\pi\mathbb{Z}+(R^2(a+b)+iR)\mathbb{Z}.$$
Since $R>0$ and $a,b \in \mathbb{R}$ are arbitrary, we see that for this connection $\varphi$, all conformal classes of a genus $1$ torus can be realized using circles in $\mathbb{R}^2$.
\\
\\ 
The appearance of what we call the \textit{inverse isoperimetric problem} in Example 2 highlights a fundamental challenge for any approach to Question~\ref{Question} on manifolds with trivial homology: the necessity of relating the holonomy angle of a Jordan curve to its length. In both \cite{Pinkall} and Example 2, the holonomy angle is expressed in terms of the enclosed area via Theorem~\ref{thm: holonomy in terms curvat}. Consequently, the isoperimetric inequality becomes the natural tool to bridge the holonomy angle and the curve's length. In the examples throughout this section, the role of the isoperimetric inequality is played by the structural relation between the contour integral of a vector field and the length of the closed curve.
\\
\\
The Ambrose-Singer theorem, which relates holonomy and curvature, could potentially be leveraged to refine the results of this section by demonstrating that the only connections for which Question~\ref{Question} is answered negatively are those with zero curvature. More specifically, in such a case, one would need to prove that for a given holonomy angle (at a point) and a given length, there exists a curve that realizes these values. 
\section{Some examples on the sphere}\label{Sec: examples sphere}
In this section, we will utilize the properties of connections on principal $S^1$-bundles and their relationship with connections on Hermitian line bundles, as outlined in Section~\ref{Sec: basic DG}.
\\
\\
Given a principal $S^1$-bundle $P$ with a connection $1$-form $\varphi \in \Omega^1(P; i\mathbb{R})$, we consider the associated Hermitian line bundle $L_P$ with its Hermitian metric, denoted $h$, and its connection, denoted $\nabla$. The $k$-fold tensor product of $L_P$ is denoted by $L_P^{\otimes k}$. This bundle is naturally endowed with an induced Hermitian metric, denoted $h^{\otimes k}$, derived from $h$. This metric, in turn, allows us to associate a principal $S^1$-bundle, which we denote $P^{\otimes k}$, to $L_P^{\otimes k}$. The connection $\varphi$ on $P$ induces a compatible connection $\nabla$ on $L_P$, which then induces a compatible connection $\nabla^{\otimes k}$ on $L_P^{\otimes k}$.
\\
\\ 
Let $\Omega$ be the curvature of the connection $\nabla$, we have
$$
\pi^*\Omega=d\varphi
$$
and that $\varphi$ is a horizontal two form, hence the previuos formula allows to recover the curvature $\Omega$ from $\varphi$.
\\
\\
If $\nabla_1$ and $\nabla_2$ are connections on vector bundles $E_1$ and $E_2$ over $M$ respectively, then the connection $\nabla$ induced by  $\nabla_1$ and $\nabla_2$ satisfies
$$
\nabla(s_1 \otimes s_2)=(\nabla_1 s_1) \otimes s_2+ s_1 \otimes (\nabla_2 s_2),
$$
where $s_i$ is a section of $E_i$. From this formula is easy to deduce that the curvature $\Omega$ of $\nabla$ is related with the curvatures $\Omega_1$ and $\Omega_2$ of $\nabla_1$ and $\nabla_2$ in the following way
\begin{equation}\label{eq:curvature}
\Omega(X,Y)(s_1 \otimes s_2)=\left(\Omega_1(X,Y)(s_1)\right) \otimes s_2+ s_1\otimes \Omega_2(X,Y)(s_2).
\end{equation}
where $X,Y$ are tangent vector of $M$ at $p$ and  $s_1$ and $s_2$ are vectors in the fiber of $E_1$ and $E_2$ at $p$.
\\
\\
If we denote $\Omega^{\otimes k}$ the curvature of $\nabla^{\otimes k}$, equation (\ref{eq:curvature}) implies that 
\begin{equation}\label{eq:Curvatura tensorial}
\Omega^{\otimes k}=k \Omega
\end{equation}
where now $\Omega$ denotes the curvature of the connection $\nabla$ on $L_P$ induced by the connection $1$--form $\varphi$ of $P$.
\\
\\
In~\cite[Proposition 1]{Pinkall}, Pinkall demonstrates that the Hopf torus associated with a smooth Jordan curve $c$ in $S^2$ is biholomorphic to $\mathbb{C}/(2\pi \mathbb{Z}+(L/2+i A/2)\mathbb{Z})$, where $L$ and $A$ are the length and signed area of the Jordan curve $c$. This result is a consequence of Theorem~\ref{thm: holonomy in terms curvat} and the symmetries of the Hopf fibration, which ensure that the signed area $A$ corresponds to the holonomy angle of the curve $c$.
\\
\\
Building on this, Lemma~\ref{prop: torus associated to connection}, Theorem~\ref{thm: holonomy in terms curvat}, and Equation~(\ref{eq:Curvatura tensorial}) imply that if we take $P:=(\pi: S^3 \to S^2)$ to be the Hopf fibration as described in~\cite{Pinkall}, then changing from $P$ to $P^{\otimes k}$ and adjusting the connections as explained above, transforms the conformal class associated with $c$ from $\mathbb{C}/(2\pi\mathbb{Z}+(A/2+iL/2)\mathbb{Z})$ to $\mathbb{C}/(2\pi\mathbb{Z}+(kA/2+iL/2)\mathbb{Z})$.
\\
\\
According to \cite[Section 3]{Pinkall} the $A$ and $L$ corresponding to signed areas and lengths of (smooth) Jordan curves of $S^2$ have as only restrictions the inequalities $0\leq A \leq 2\pi$ and the isoperimetric inequality 
$$
L^2-4\pi A-A^2 \geq 0.
$$  
According to \cite[Section 3]{Pinkall}, the signed areas ($A$) and lengths ($L$) of smooth Jordan curves on $S^2$ are constrained only by $0 \leq A \leq 2\pi$ and the isoperimetric inequality:
$$L^2 - 4\pi A - A^2 \geq 0.$$
The homothety $(A/2, L/2) \mapsto (k A/2, L/2)$ transforms the region defined by these inequalities into a new region. When this transformed region is joined with its reflection across the y-axis, it can be observed to contain a fundamental domain of the $SL(2,\mathbb{Z})$ action on the upper half-plane of the complex plane, where the coordinates are given by $A/2 + iL/2$. This leads to the following proposition, which answers Question~\ref{Question} positively.
\begin{prop}
Let $P:=(\pi:S^3 \to S^2)$ be the Hopf fibration, endowed with the connection $1$--form $\varphi$ inherited from the Riemannian metric of $S^3$. Then, every conformal class of a compact surface of genus $1$ can be obtained as the pullback of the principal $S^1$--bundle $P^{\otimes k}$ endowed with the connection $\varphi^{\otimes k}$ to a smooth Jordan curve $c$ in $S^2$, using Lemma~\ref{prop: torus associated to connection}. Moreover, if the Jordan curve $c$ has signed area $A$ and length $L$, the complex torus $c^*(P^{\otimes k})$ is biholomorphic to $\mathbb{C}/(2\pi\mathbb{Z}+(L/2+i k A/2)\mathbb{Z})$.
\end{prop}
\begin{rem}\label{rem: sphere}
The results of this section illustrate the exploration of the moduli space $\mathscr{M}$ through the variation of the bundle and connection over a fixed base geometry. For a fixed Jordan curve $c \subset S^2$, its length $L$ provides one stable parameter in $\mathscr{M}$. The second parameter is the holonomy angle, which varies according to the choice of the $S^1$--bundle and its connection. Despite these variations, the holonomy remains a rigid geometric quantity: for any given connection on a fixed line bundle, the holonomy is uniquely determined by the enclosed area of the curve. Consequently, by transitioning through different bundle structures, we effectively vary the holonomy parameter of the moduli space while the geometry of the base curve $c$ remains fixed.
\end{rem}
\section{Conclusions and perspectives}\label{Sec: Conclusions and P}
In this article, we have presented a natural generalization of Garsia's conjecture as formulated in Question~\ref{Question}. By examining this conjecture through the lens of trivial homology cases ($H_1(M, \mathbb{Z}) = 0$), we have identified a fundamental topological rigidity: the holonomy of a Jordan curve is uniquely determined by the curvature flux through any bounding surface. This insight allows us to interpret the parameters of the moduli space $\mathscr{M}$ in terms of tangible geometric quantities—specifically length and holonomy—as illustrated in Remarks~\ref{rem:Parameter length and area}, \ref{rem:holonomy as an integral of curvature}, \ref{rem:parametrization of moduli Rn}, and \ref{rem: sphere}. 
\\
\\
We propose that Lemma~\ref{prop: torus associated to connection} provides a robust framework for constructing complex tori within a differential geometric context, offering a path to study elliptic curves and their moduli without a purely algebraic approach. 
\\
\\
The present work establishes the "simply connected" theory of the conjecture, where the area-holonomy relationship is rigid. A natural and more challenging direction for future research involves manifolds with non-trivial homology, where the holonomy is no longer uniquely determined by the enclosed area. In a work in preparation, we shall investigate how the holonomy varies for curves with non-trivial homology, a study which we expect will shed further light on the broader generalization of Garsia's conjecture. Furthermore, it would be of great interest to explore how this methodology can be extended to higher genus surfaces and higher-dimensional manifolds.\bibliographystyle{alpha}
\bibliography{literatur}
\end{document}